\documentclass[12pt]{amsart}
\usepackage{amsfonts, amssymb,amsmath, amsthm, amsxtra, latexsym, mathrsfs, amscd}
\usepackage{amsmath, amssymb, amsthm, times, float}
\usepackage{graphics}
\usepackage[latin1]{inputenc} %

\newtheorem{theo+}{Theorem}[section]
\newtheorem{prop+}[theo+]{Proposition}
\newtheorem{coro+}[theo+]{Corollary}
\newtheorem{lemm+} [theo+]{Lemma}
\newtheorem{deep+}  [theo+]  {Deep Result}
\newtheorem{fact+}  [theo+]  {Fact}
\theoremstyle{definition}
\newtheorem{exam+}  [theo+]  {Example}
\newtheorem{rema+}  [theo+]  {Remark}
\newtheorem{defi+}  [theo+]  {Definition}
\newtheorem{xca+}[theo+]{Exercise}

\numberwithin{equation}{section}

\usepackage{color}

\def\draft{\centerline{(Draft {\the \day}/{\the\month} \the \year.)}}
\bigskip

\def\refn#1.#2{\expandafter\def\csname#1\endcsname{[#2]}}
\def\refnr#1.{\csname#1\endcsname}

\def\fa{\mathfrak a}

\def\fg{\mathfrak g}
\def\fk{\mathfrak k}
\def\fh{\mathfrak h}

\def\fm{\mathfrak m}
\def\fn{\mathfrak n}
\def\fp{\mathfrak p}

\def\a{\alpha}

\def\Claminv2{|C(\Lambda)|^{-2}}

\def\de{d\varepsilon}

\def\Aa2D{A^{\a,2}(D)}
\def\bAa2D{\overline{A^{\a,2}(D)}}
\def\Ab2D{A^{\beta,2}(D)}
\def\bAb2D{\overline{A^{\beta,2}(D)}}

\def\Norm#1_#2{\Vert#1\Vert_{#2}}

\def\phipl12{\phi_{p_{l_1}, p_{l_2}}}
\def\phip01{\phi_{p_{0}, p_{0}}}

\def\a{\alpha}

\def\Claminv2{|C(\Lambda)|^{-2}}

\def\sig{\sigma}

\def\de{d\varepsilon}

\def\Aa2D{A^{\a,2}(D)}
\def\bAa2D{\overline{A^{\a,2}(D)}}
\def\Ab2D{A^{\beta,2}(D)}
\def\bAb2D{\overline{A^{\beta,2}(D)}}

\def\phipl12{\phi_{p_{l_1}, p_{l_2}}}
\def\phip01{\phi_{p_{0}, p_{0}}}

\def\bc{\mathbb C}
\def\br{\mathbb R}

\def\bh{\mathbb H}

\def\alg/{algebra}

\def\Alg/{Algebra} 

\def\alt/{alternative} 
\def\anal/{analytic}
\def\analfunc/{\anal/\ \func/}
\def\Ans/{\it Answer. \normal}
\def\ass/{associative}
\def\nass/{non-\ass/}
\def\autom/{automorphism}
\def\homom/{homomorphism}
\def\isom/{isomorphism}
\def\bdd/{bounded}
\def\Bdd/{Bounded}
\def\bddsymdom/{bounded \sym/ \dom/}
\def\Cartdom/{Cartan \dom/}
\def\bdry/{boundary}
\def\bsd/{\bdd/ \symdom/}
\def\bv/{boundary value}
\def\cf/{{\it cf}\.}
\def\Cf/{{\it Cf}\.}
\def\charr/{character}
\def\coeff/{coefficient}
\def\comm/{commutative}
\def\cpct/{compact}
\def\compl/{complex}
\def\comp/{complex}
\def\Comp/{Complex}
\def\conf/{conformal}
\def\conj/{conjugate}
\def\conn/{connect}
\def\cont/{continuous}
\def\conv/{converge} 
\def\convc/{convergence}
\def\convt/{convergent}
\def\convx/{convex}
\def\coord/{coordinate}
\def\lcoord/{local coordinate}
\def\Corr/{Corresponding}
\def\corr/{corresponding}
\def\corrd/{correspond}
\def\cov/{covariant}
\def\decomp/{decomposition}
\def\deco/{decompose}
\def\diff/{different} 
\def\Diff/{Different} 
\def\dimn/{dimension} 
\def\distr/{distribution} 
\def\div/{diverge} 
\def\dom/{domain}
\def\eg/{\hbox{\it e.g}\.}
\def\eigenf/{eigen\-\func/}
\def\eigensp/{eigen\-space}
\def\eigenv/{eigen\-value}
\def\eq/{equation}
\def\equa/{equation}
\def\de/{\diff/ial \equa/}
\def\do/{\diff/ial operator}
\def\ode/{ordinary \de/}
\def\pde/{partial \de/}
\def\pdo/{partial \diff/ial operator}
\def\psdo/{pseudo \diff/ial operator}
\def\fin/{finite}
\def\Ex/{\it Example.\ \normal}
\def\Exnr#1/{\it Example #1.\ \normal}
\def\foll/{follow}
\def\follg/{following}
\def\Follg/{Following}
\def\func/{function}
\def\Func/{Function}
\def\Fonc/{Fonc\-tion}
\def\fonc/{fonc\-tion}
\def\Funk/{Funk\-tion}
\def\funk/{Funk\-tion}
\def\gen/{general}
\def\har/{harmonic}
\def\Hint/{\it Hint. \normal}
\def\hist/{historic}
\def\histcl/{historical}
\def\hol/{holo\-morphic}
\def\homog/{ho\-mo\-ge\-ne\-ous}
\def\hyp/{hyper\-bolic}
\def\hyperg/{hyper\-geometric}
\def\ie/{\hbox{\it i.e.}}
\def\iff/{if and only if}
\def\ineq/{inequality}
\def\infra/{{\it inf\-ra}}
\def\ultra/{{\it ult\-ra}}
\def\Inpart/{In particular}
\def\inpart/{in particular}
\def\instof/{instead of}
\def\interps/{interpolation space}
\def\interp/{interpolation}
\def\Interp/{Interpolation}
\def\interpr/{Interpretation}
\def\Intr/{Introduction}
\def\intv/{interval}
\def\inv/{invariant}
\def\invc/{invariance}
\def\Iowords/{In other words}
\def\iowords/{in other words}
\def\ipr/{inner product}
\def\irred/{irreducible}
\def\lb/{line bundle}
\def\lin/{linear}
\def\lhs/{left hand side}
\def\rhs/{right hand side}
\def\loc/{local}
\def\math/{mathematic} 
\def\mathcn/{\math/ian}
\def\manif/{manifold}
\def\meas/{measure}
\def\measl/{measurable}
\def\mero/{mero\-morphic}
\def\mon/{monomial}
\def\monog/{monogenic}
\def\mult/{multiple}
\def\multy/{multiply}
\def\multn/{multiplication}
\def\nas/{necessary and sufficient}
\def\nbd/{neighborhood}
\def\neg/{negative}
\def\nondeg/{nondegenerate}
\def\Oohand/{On the other hand}
\def\oohand/{on the other hand}
\def\Oonhand/{On the one hand}
\def\oonhand/{on the one hand}
\def\oper/{operator}
\def\orth/{ortho\-gonal}
\def\orthon/{ortho\-normal}
\def\otoh/{on the other hand}
\def\quat/{quaternion}
\def\pp/{\hbox{a. e.}}
\def\psh/{plurisubharmonic}
\def\pol/{polynomial}
\def\pot/{potential}
\def\pos/{positive}
\def\princ/{principle}
\def\prob/{probability}
\def\proj/{projective}
\def\projn/{projection}
\def\Proof/{\it Proof:\normal}
\def\Rem/{\it Remark\normal}
\def\Remnr#1/{\it Remark\ \normal #1. }
\def\rep/{representation}
\def\reps/{representations}
\def\meta/{metaplectic representation}
\def\repr/{reproducing}
\def\reprker/{reproducing kernel}
\def\resp/{respective} 
\def\resply/{respectively}
\def\restr/{restriction}
\def\sa/{self-adjoint}
\def\st/{such that}

\def\sol/{solution}
\def\ru/{space}
\def\sph/{spherical}
\def\ssp/{sub\ru/}
\def\sym/{symmetric}
\def\Sym/{Symmetric}
\def\symb/{symbol}
\def\symbc/{symbolic}
\def\symdom/{\sym/ domain}
\def\symp/{symplectic}
\def\Theor#1/{\fet Theorem #1.\ \normal}
\def\Lem#1/{\fet Lemma #1.\ \normal}
\def\Lemma/{\fet Lemma.\ \normal}
\def\topl/{topology}
\def\topll/{topological}
\def\transf/{transform}
\def\transl/{translation}
\def\transfn/{transformation}
\def\transv/{transvectant}
\def\trig/{trigonometric}
\def\tril/{trilinear}
\def\trilf/{trilinear form}
\def\uhp/{upper halfplane}
\def\uhs/{upper halfspace}
\def\vb/{vector bundle}
\def\vf/{vector field}
\def\vsp/{vector space}
\def\wrt/{with respect to}
\def\Wlog/{Without loss of generality}
\def\a{\alpha}

\def\sig{\sigma}

\def\Ab/{Abel}
\def\Ban/{Banach}
\def\Bansp/{\Ban/ space}
\def\Belt/{Bel\-tra\-mi}
\def\Berg/{Berg\-man}
\def\Bern/{Ber\-nou\-lli}
\def\Berz/{Berezin}
\def\Bess/{Bessel}
\def\Cart/{Car\-tan}
\def\Cay/{Cay\-ley}
\def\CG/{Clebsch-Gordan}
\def\Cl/{Clifford}
\def\CR/{Cauchy-Rie\-mann}
\def\Dir/{Dirichlet}
\def\Eucl/{Euclide}
\def\Eucln/{Euclidean}
\def\F/{Fourier}
\def\Hank/{Hankel}
\def\Hankf/{\Hank/ form}
\def\Herm/{Hermite}
\def\Hilb/{Hilbert}
\def\Hilbs/{Hilbert space}
\def\Hilbsp/{Hilbert space}
\def\HS/{Hilbert-Schmidt}
\def\Lag/{La\-grange}
\def\Lap/{La\-place}
\def\LapBelt/{\Lap/-\Belt/}
\def\Leb/{Lebesgue}
\def\Marc/{Mar\-cin\-kie\-wicz}
\def\Moeb/{Moebius}
\def\Moebt/{Moebius transformation}
\def\Moebtransfn/{Moebius transformation}
\def\Pla/{Plan\-che\-rel}
\def\Poin/{Poin\-car\'e}
\def\Riem/{Rie\-mann}
\def\Riemn/{\Riem/ian}
\def\psRiemn/{pseudo-\Riem/ian}
\def\Riems/{Rie\-mann surface}
\def\Schroe/{Schr\"odinger}
\def\Weier/{Weier\-strass}
\def\anal/{analytic}
\def\bsd/{bounded symmetric domain  }
\def\bdd/{bounded}
\def\calc/{calculation}\def\conj{conjugate}
\def\calci/{calculating}\def\eg{e.g.}
\def\conj/{conjugate}
\def\deco/{decomposition}
\def\eg/{e.g.}
\def\fct/{function}
\def\gp/{group}
\def\hw/{highest weight}
\def\hwv/{highest weight vector}
\def\hwvs/{highest weight vectors}
\def\lw/{lowest weight}
\def\lwv/{lowest weight vector}
\def\lwvs/{lowest weight vectors}
\def\hds/{holomorphic discrete series}
\def\iff/{if and only if}
\def\inv/{invariant}
\def\irrde/{irreducible decomposition}
\def\meas/{measure}
\def\transf/{transform}
\def\rep/{representation}
\def\resp/{respectively}
\def\inters/{intertwines}
\def\interg/{intertwining}
\def\meta/{metaplectic representation}
\def\qu/{quaternion}
\def\rep/{representation}
\def\symdom/{ symmetric domain}
\def\st/{such that}
\def\shd/{subhead}
\def\transf/{transform}
\def\wrt/{with respect to}

\def\Norm#1#2#3{\Vert#1\Vert^{#3}_{{#2}}}

\begin{document}
\def\abstractname{Abstract}
\def\chrefname{References}

\title[Restriction of unitary representations
]{
Discrete components in restriction 
of unitary representations of rank one 
semisimple Lie groups
 }
\author{ Genkai Zhang}
\address{Mathematical Sciences, Chalmers University of Technology and
Mathematical Sciences, G\"oteborg University, SE-412 96 G\"oteborg, Sweden}
\email{genkai@chalmers.se}
\thanks{Research partially supported by 
the Swedish
Science Council (VR)}
\begin{abstract}
We consider spherical principal series
representations 
of the semisimple Lie group
of rank one 
$G=SO(n, 1; \mathbb K)$,
$\mathbb K=\br, \bc, \bh$.
 There is a family of
 unitarizable representations $\pi_{\nu}$ of $G$
for $\nu$ in an interval
on $\mathbb R^+$, the so-called complementary series,
and subquotients or subrepresentations of $G$
for $\nu$ being negative integers.
We consider the restriction of 
$(\pi_{\nu}, G)$ under the subgroup
$H=SO(n-1, 1; \mathbb K)$.
We prove the appearing of discrete
components.
The corresponding  results
for the exceptional Lie group $F_{4(-20)}$ and its
 subgroup $Spin(8,1)$
are also obtained.
\end{abstract}

\maketitle

\baselineskip 1.35pc

\section{Introduction}

The study  of direct components
in the restriction to a subgroup $H\subset G$
of a representation $(\pi,  G)$ is one of major
subjects in representation theory. 
Among  representations of  a semisimple Lie
group $G$ there are two somewhat opposite
classes, the discrete series
and the complementary series;
the former appear in the decomposition of $L^2(G)$
 and can be treated 
algebraically, whereas the latter 
do not contribute to the decomposition
and their study involves  more analytic issues.
The study of restriction of discrete
series representations has been studied
intensively; see e.g. \cite{Ors-Var, kob-rest-2005}
and references therein.
Motivated by some related questions 
of \cite{Bergeron-imrn, Cloz-Berg}
 Speh and
Venkataramana \cite{Speh-Venk-2} 
studied the restriction
of a complementary series
representation of $SO(n, 1)$
under the subgroup $SO(n-1, 1)$.
It is approved there, 
 for relatively small parameter $\nu$ (in our parametrization),
 the complementary series $\pi_\nu$
of  $SO(n-1, 1)$
appears discretely 
in the complementary series
 $\pi_\nu$ of  $SO(n, 1)$ 
with the same
parameter $\nu$.
They construct the imbedding
of the complementary series
of  $SO(n-1, 1)$  into 
 $\pi_\nu$ of  $SO(n, 1)$ 
by  using  non-compact
realizations of the representations 
as spaces of distributions on Euclidean spaces and 
by extending distributions on $\mathbb R^{n-2}$
to $\mathbb R^{n-1}$. 
Similar 
results are also obtained for 
complementary
series of differential
forms.
For $n=3$
the same result
is proved \cite{Speh-Venk-1} 
by using the compact picture
and $SU(2)$-computations;
see also \cite{Mukunda} where
a full decomposition is found.

In the present paper we shall study
the branching of complementary
series of $G$ for 
all rank one Lie groups $G$
under a symmetric subgroup $H$. More precisely
we prove the appearance
of discrete components
for $G=SO(n, 1; \mathbb K)$,  $H=SO(n-1, 1; \mathbb K)$,
with $\mathbb F=\br, \bc, \bh$
being the fields of real, complex, quaternion
numbers, or for $G=F_{4(-20)}$
 and $H=Spin(8, 1)\subset G$.
We shall use the compact realization
of the spherical principal series $\pi_\nu$
on the sphere $S=K/M$ in $\mathbb F^n$.
We prove that for appropriate
small parameter  $\nu$
the natural restriction map
of functions   on $S$ in $\pi_\nu$
to the lower dimensional sphere $S^\flat$
in $\mathbb F^{n-1}$ 
defines
a bounded operator onto a complementary series
$\pi_\nu^\flat$ of $H$.
The proof requires
rather detailed study of the
restriction to $S^{\flat}\subset S$  of spherical
harmonics on $S$.

The \reps/ $\pi_{\nu}$ for certain integers $\nu$
have also
 unitarizable subquotients or subrepresentations. Some of them 
are discrete series
for $SU(2, 1)$.
We shall find  irreducible
components for the \reps/
under the subgroup $H$.
One easiest case
is the subrepresentation $\pi_0^{\pm}$ 
(or  $\pi_{2n+2}^{\pm}$ as quotient)
of the group $SU(n, 1)$. The space
 $\pi_0^{\pm}$ consists
 of holomorphic respectively antiholomorphic
polynomials
on $\mathbb C^n$ modulo  constant functions.
It  can also
be treated by using the analytic
continuation of scalar holomorphic
discrete series at the reducible point
\cite{FK},
and some general decomposition results
have been obtained in \cite{kobayashi-06}.
These representations
are also of special interests in 
automorphic representation theory
\cite{Sarnak-claylecture}.

The main results in this paper is summarized in the following;
see Theorems 3.6, 3.9 and 4.4 below for the parameterization
of the complementary series and
precise statements.

\begin{theo+} 
Let $(G, H)$ be the pair as above,
$G=SO(n, 1; \mathbb K)$,  $H=SO(n-1, 1; \mathbb K)$
for  $\mathbb F=\br, \bc, \bh$,
 or $G=F_{4(-20)}$, and  $H=Spin(8, 1)\subset G$.
Let $\rho_G$ and $\rho_H$ be the
corresponding half sum of positive roots. Suppose
 $(\pi_{\nu}, G)$ is a complementary
series representation of $G$. We can assume
up to Weyl group symmetry that $\nu<\rho_G$.
\begin{enumerate}
\item
The restriction of $(\pi_{\nu}, G)$
on $H$ contains a discrete component 
$(\pi_{\mu}^\flat, H)$
if $\nu <\rho_H$, and $\mu=\nu$ in our parameterization.
\item
 Consider the analytical continuations
of  $(\pi_{\nu}, \fg)$
$(\pi_{\nu}^\flat, \fh)$,
 and the corresponding unitarizable
quotients $(\mathcal W_k, \pi_{\nu(k)}, \fg)$ and
$(\mathcal V_k, \pi_{\mu(k)}, \fg)$ for $\nu=-k$,
$\mathbb F=\mathbb R$ 
or $\nu=-2k$, $\mathbb F=\mathbb C, \mathbb H$.
Then the restriction  of 
$(\mathcal W_k, \pi_{\nu(k)}, \fg)$ 
to $\fh$ contains the
discrete component 
 $(\mathcal V_k, \pi_{\mu}, \fh)$
with the same parameters $\mu=\nu(k)$.
\end{enumerate}
\end{theo+}

We note that our result can be understood heuristically 
as certain  boundedness property 
of the restriction map from certain Sobolev spaces
 on $S$ to those on $S^\flat$.
Indeed for
small parameter $\nu$ 
the space $\pi_{\nu}$ consists of distributions
on $S$ whose fractional  differentiations
are in $L^2(S)$, i.e.,
they are functions with certain smooth conditions. It is thus
expected that their restriction
on the subsphere  $S^\flat$
would make sense in proper Sobolev spaces. (However
the precise space of complementary series
is not the usual Sobolev space; only $L^2$-conditions
for the differentiations of functions are required.)
A precise formulation can be done and we hope to return to it in future.
 We remark also that the
study of the norm estimates
of the restriction of the spherical harmonics
on lower dimensional
spheres can be put
into a general context as the study
of growth of $L^p$-norm
of restriction on totally geodesic
submanifolds of eigenstates of Laplace-Beltrami
operators on Riemannian manifold; see  \cite{Burg-etal}.
Our
results here give  precise estimates
of the  $L^2$-norm
of the restriction. They might have independent
interests on their own from a view point
of harmonic analysis.
They might also shed light on
the  study of $L^p-L^q$ estimate 
of the above restriction problem
for general compact manifolds.

I would like to thank 
 B. Speh
and T. N. Venkataramana 
for  some 
 correspondences
and stimulating
discussions during an AIM work ``Branching
problems in unitary representations'',
MPIM, Bonn, July 2011; 
they have obtained  similar results
for the case $SO(n, 1)$ and $SU(n, 1)$
earlier in an unpublished manuscript.
I thank also T. Kobayashi, 
R. Howe, B. \O{}rsted and J. Vargas for
several conversations. I am  grateful
for the anonymous referee for expert
criticisms and for suggesting
many improvements.

\section{Preliminaries}
\subsection{Classical rank one groups}
Let $\mathbb F=\mathbb R, \mathbb C, \mathbb H$
be the real, complex and quaternionic numbers.
Denote $G:=SO_0(n, 1; \mathbb F)=SO_0(n, 1),
SU(n, 1), Sp(n, 1)$
the connected component
of the group $GL(\mathbb F^{n+1})$
 of $\mathbb F$-linear transformations
on $\mathbb F^{n+1}$ preserving
the quadratic form $|x_1|^2 +\cdots +|x_n|^2 -|x_{n+1}|^2$,
with $\mathbb F$ acting on the right.
The group $K:=SO_0(n), S(U(n)\times U(1)),
Sp(n)\times Sp(1)$ is a maximal
compact subgroup of $G$ and $G/K$
is a Riemannian symmetric space of rank one
which can further be realized as the unit ball
in $\mathbb F^n$. 
Elements in $G$ and $\fg$ will be written
as $(n+1)\times (n+1)$ block $\mathbb F$-matrices 
$$
\begin{bmatrix} a&b\\
c&d
\end{bmatrix}
$$
where $a, b, c, d$ are of size $n\times n,
n\times 1, 1\times n, 1\times 1$, respectively.

Let $\fg=\fk +\fp$ be the corresponding
Cartan decomposition.
We fix
$$
H_0=\begin{bmatrix} 0&e_1\\
e_1^T&0
\end{bmatrix}
$$
in $\fp$, where $e_1$
is the standard basis vector
and  $e_1^T$ its transpose, and let $\fa=\mathbb RH_0\subset \fp$. 
Then $\fa$
 is a maximal abelian subspace of 
$\fp$. The root space decomposition of $\fg$ under  $H_0$ 
is 
$$
\fg=\fg_{-1} +(\fa+ \fm) +\fg_{1} 
$$
with roots $\pm 1, 0$ 
if $\mathbb F=\mathbb R$, and 
$$
\fg=\fg_{-2} + \fg_{-1} +
(\fa+ \fm) + \fg_{1} +\fg_{2} 
$$
with roots $\pm 2, \pm 1, 0$,
if $\mathbb F=\mathbb C, \bh$.  Here $
\fm\subset \fk$ is the zero root space
in $\fk$.
We denote
$$\fn=\fg_{1}, \quad\fn=\fg_{1}
+\fg_2$$ 
the sum of the positive root spaces,
in the respective cases.
 Thus 
$\fm+\fa +\fn$
is a maximal parabolic subalgebra
of $\fg$.
Let $\rho$
be the half sum of positive roots. We identify
$\fa^\ast$ with $\mathbb C$ via
$\lambda\to \lambda(H_0)$ and  we write $\rho(H_0)=\rho$.
We have
\begin{equation}
  \label{eq:realrho}
\rho=
\begin{cases}
\frac{n-1}2, \quad \mathbb F=\br\\
n, \quad \mathbb F=\bc\\
2n+1, \quad \mathbb F=\bh.
\end{cases}  
\end{equation}
Denote $M, A, N$ 
the corresponding 
subgroups with Lie algebras
$\fm, \fa, \fn$. Then  $M=SO_0(n-1)=SO(n-1)$, $\tilde{SU(n-1)}$
a double cover of ${SU(n-1)}$, or $Sp(n-1)\times Sp(1)$,
and $MAN$ is a maximal parabolic
subgroup of $G$.

\subsection
{Decomposition of $L^2(K/M)$} 
We identify $\fp$ with $\mathbb F^n$
and  normalize
the $K$-invariant inner product on $\fp$ so that
$H_0$ is a unit vector.
The homogeneous space $K/M$ is then
the unit sphere $S:=S^{dn-1}$ in 
$\fp=\mathbb F^n$ with $M$
being the isotropic subgroup of the base point $H_0\in \frak \fp$,
with $d=\text{dim}_{\mathbb R}F=1, 2, 4.$
We denote $dx$ the normalized area measure
on $S$ and $L^2(S)$ the corresponding $L^2$-space.
For $n=1$ the decomposition of $L^2(K/M)$ is
well-known and elementary, so we assume $n>1$. 
Let $W^p$ be the space of
spherical harmonics on $S$. For $\mathbb F=\bc$
let  $W^{p, q}$ be the spherical harmonics
of degree $p+q$ on $\mathbb C^n$ and holomorphic of degree $p$
and antiholomorphic of degree $q$.
If  $\mathbb F=\bh$, then $K=Sp(n)\times Sp(1)$, 
and its representations are of the form $\tau_1\times \tau_1$,
which will be written as $(\tau_1, \tau_2)$
and further identified with their highest
weights.
The root system of $Sp(n)$
is of type $C$ and let $\alpha_1, \cdots, \alpha_{n-1}, \alpha_n$
be the simple roots with $\alpha_n$ the longest one.
Denote $\lambda_1, \cdots, \lambda_n$ the corresponding
fundamental weights with $\lambda_1$
the defining representation on $\mathbb C^{2n}$.
 For $Sp(1)=SU(2)$ the representation
on symmetric tensor power $\odot^q(\mathbb C^2)=\mathbb C^{q+1}$
will be written just as $q$
 for simplicity.
Denote $W^{p, q}$
 the representation $(q\lambda_1 +\frac{p-q}2\lambda_2, q)$
of $K=Sp(n)\times Sp(1)$.

Recall \cite{Kostant-BAMS, Johnson-Wallach}
\begin{equation}
L^2(S)=\sum_{\tau}^\oplus W^\tau, \quad
W^\tau=\begin{cases}
W^p, \, p\ge 0\, & \mathbb F=\br\\
W^{p, q}, \, p, q\ge 0\, &\mathbb F=\bc\\
W^{p, q}, \, {p\ge q\ge 0, p-q \,\text{even}}
 \, &\mathbb F=\bh
\end{cases}
\end{equation}
Here and in the following we denote a general representation of $K$
 by $\tau$.
 The subspace 
 $(W^\tau)^M$ of  $M$-fixed vectors
is one dimensional 
$$
(W^\tau)^M =\mathbb C \phi_{\tau}
$$
where $\phi_{\tau}$ is normalized by $\phi_{\tau}(H_0)=1$.
They depend only on the first variable $x_1\in \bh$
of $x=(x_1, \cdots, x_n)$,  and will also be written
 as $\phi_{\tau}^n(x_1)$.
We recall some explicit formulas for
them obtained in  \cite[Theorem 3.1]{Johnson-Wallach}. 
 (Note that in the formula for $\psi_{p, q}$
and $e_{p, q}$ in  \cite[p.144-147]{Johnson-Wallach}
the term $\frac {-p-q}2$ should be
$\frac {-p+q}2$.) Those polynomials
are obtained
as polynomial solutions to differential
equations.
A variant of these
polynomials will be constructed in Lemma 3.3.
\begin{lemm+} The polynomials  $\phi_{\tau}^n$
are given as follows:
\begin{enumerate}
\item $\mathbb F=\mathbb R$,  $x_1=\cos\xi$,
\begin
{equation*}
\phi_p^{n}(x_1):=
\cos^p \xi F(-\frac p2, 
-\frac {p-1}2, \frac{n-1}2, -\tan^2\xi);
\end{equation*}
\item $\mathbb F=\mathbb C$,
$x_1=e^{i\theta }\cos\xi$,
\begin{equation*}
\phi_{p, q}^{n}
(x_1)=
e^{i\theta(p-q)}\cos^{p+q} \xi
 F(-p, 
-q, n-1, -\tan^2\xi);
\end{equation*}
\item
$\mathbb F=\mathbb H$,   $x_1=\cos\xi e^{\theta y}
=\cos\xi
 (\cos \theta + y\sin \theta)$
 in  quaternionic  polar coordinates,
$y\in H$ being purely imaginary (i.e. in $\mathbb Ri+\mathbb Rj+\mathbb Rk$)
and $|y|=1$,
\begin{equation*}
\phi_{p, q}^{n}(x)=
\phi_{p, q}^{n}(x_1):=
\frac{\sin(q+1)t}
{\sin t}
\cos^{p} \xi\,
 F(-\frac {p-q}2, 
-\frac {p+q+2}2 , 2(n-1), -\tan^2\xi).
\end{equation*}
\end{enumerate}
\end{lemm+}

Here $F(a, b, c, x)$
is the Gauss hypergeometric function ${}_2F_1$,
 $$
F(a, b, c, x)=\sum_{m=0}^\infty
\frac{(a)_m(b)_m}{(c)_m}\frac {x^m}{m!}
$$
 and $(a)_m=\prod_{j=0}^{m-1}
(a+j)$ is the Pochammer symbol. Note that all $\phi$-functions
above are Jacobi polynomials \cite{Szeg} in $t=2|x_1|^2-1$ in the interval
$(-1, 1)$.

We put the upper-index the dimension $n$ as we shall also
treat it as a variable.

In particular we have, by Schur's orthogonality relation,
\begin{equation}
  \label{eq:schur-1}
\Vert \phi_\tau
\Vert^2=\frac 1{\text{dim}(W^\tau)}.
\end{equation}
$\text{dim}(W^\tau)$ can be evaluated by the Weyl's dimension
formula: Let $\{\alpha\}$ be the root system of $\fk$
with $\{\alpha>0\}$ the positive roots
and $\rho_{\fk}
$ the half sum of the positive roots,
$$
\text{dim}(W^\tau)=\prod_{\alpha>0}
\frac{
\langle\tau+\rho_{\fk}, \alpha\rangle
}
{\langle\rho_{\fk}, \alpha\rangle
}
;
$$
see e.g. \cite{He2}.

We shall also need a  general integral formula: If $f(x)=g(y)
h(z)$, $x=(y, z)$ are functions
on $\mathbb R^{m}$ with separated variables  $y\in \mathbb R^k$
and $z\in \mathbb R^{m-k}$ with $dy$ the Lebesgue 
measure then we have
\begin{equation}
  \label{eq:polar}
\begin{split}
\int_{S^{m-1}}
f(x) dx
&=
\frac{2\Gamma(\frac m2)}
{
\Gamma(
\frac{k}2
)
\Gamma
(
\frac{m-k}2)
\omega_{k-1}
}
\int_{B^k} g(y) (1-|y|^2)^{\frac 12
 (m-k-2)}
\\
&\quad\times
\left(\int_{S^{m-k-1}} h( (1-|y|^2)^{\frac 12}z)  dz\right)
 dy
\end{split}
\end{equation}
where $B^k$ is the unit ball in $\mathbb R^k$,
$dx$ and $dz$ are the area measures 
on the respective spheres normalized  
with total areas being  $1$,
and
$\omega_{k-1}=\frac{2\sqrt{\pi}^k}{\Gamma(\frac k2)}$
 is the Lebesgue area of the 
sphere in $\mathbb R^k$ (we shall  need $k=1, 2, 4$ only);
 see e.g. \cite[1.4.4 (1)]{Rudin-ball}
for the case of even $m$ and $k$. 
Thus the square norm $\Vert \phi_\tau\Vert^2$ can also be proved
by using the known integral formulas
for Jacobi  polynomials.
 However
we shall use mostly the Weyl's dimension formula 
whenever possible as it is conceptually
clearer and as their asymptotic are well-understood.

\subsection{Exceptional group $F_{4(-20)}$} Let $G$ be the 
connected Lie group
of type $F_{4(-20)}$ with Lie algebra
$\fg$.  This group has been well-studied
\cite{Johnson-2, Takahashi}.
The maximal
compact subgroup $K$ is $Spin(9)$
and the symmetric space $G/K$
can be realized 
as the unit ball in $\mathbb O^2$ with 
$\mathbb O$ being the Cayley division (octonian) algebra.
Let $\fg=\fp+\fk$ be the Cartan decomposition.
The space $\fp$ will be identified with 
$\mathbb O^2$ 
 with $\fk=spin(9)$
acting on $\mathbb O^2$ via the Spin representation.
We fix $H_0\in \mathbb O^2= \fp$
so that the positive
eigenvalues of $\text{ad}(H_0)$ in $\fg$ are
$2,  1$.
 The corresponding
multiplicities are then $7$ and $8$.
The half sum of positive roots is  $\rho=11$.
Let $\fm$ be the zero root space of $H_0$
in $\fk$, and $\fm +\fa +\fn$
the maximal parabolic subalgebra.

The algebra $\fm\subset \fk$ is $\mathfrak{spin}(7)$. 
Let $M=Spin(7)$ be the corresponding
simply connected  subgroup with Lie algebra $\fm$. 
Fix the $K$-invariant
inner product on $\fp=\mathbb O^2$ with $H_0$ being unit vector.
The homogeneous
space $K/M$ is the unit sphere $S=S^{15}$
in $\mathbb O^2=\mathbb R^{16}$. To describe
the decomposition of $L^2(S)$ under $K$ we
observe first that the space $\fp=\mathbb O^2$
is decomposed under $M$
as 
\begin{equation}
  \label{eq:r16-m}
  \fp=\mathbb O \oplus \mathbb O
=(\mathbb RH_0 \oplus\mathbb R^7)
\oplus \mathbb O
\end{equation}
with $\mathbb R^7$ being the defining representation
of $SO(7)$ and thus of $M$ via the double covering
$M=Spin(7)\to SO(7)$, and $\mathbb O$
the Spin representation of $M$.
The Dynkin diagram of $Spin(9)$
is
  \begin{center}
\begin{figure}[h]
\includegraphics{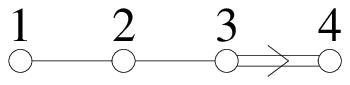}
\end{figure}
  \end{center}
with the simple roots $\alpha_1,
\alpha_2, \alpha_3, \alpha_4$. Let 
$\lambda_1, \lambda_2, \lambda_3, \lambda_4$ be the corresponding fundamental weights.
Let $
W^{p, q}
$ be the  representation
of $K$ with highest weight $\frac{p-q}2\lambda_1 + q\lambda_4$.
Then it follows \cite{Kostant-BAMS, Johnson-2}
that
\begin{equation}
  \label{eq:L-2-cayley}
L^2(S)=\oplus_{p\ge q\ge 0, \, p-q\ge 0 \, \text{even}}
W^{p, q},  
\end{equation}
and each space $W^{p, q}$ has a 
unique $M$-fixed vector $\phi_{p, q}$, $(W^{p, q})^M=\mathbb C \phi_{p, q}$,
such that $\phi_{p, q}(H_0)=1$.
To describe $\phi_{p, q}$
write elements in $\mathbb O^2$
as $x=(x_0, x_1, x_2)$
under the decomposition (\ref{eq:r16-m}), 
and write
their (partial) polar coordinates as
$r=|x|$, $\sqrt{x_0^2+\Vert x_1\Vert^2}=r\cos\xi $,
$x_0=r\cos\xi \cos\eta$
with $0\le \xi\le \frac \pi 2$, 
$0\le \eta\le  \pi $. Then 
\begin{equation*}
\begin{split}
&\quad\,\phi_{p, q}(x)=\phi_{p, q}(x_0)\\
&=
\cos^q \eta \, F(-\frac q 2, 
-\frac {q-1} 2, \frac 72; -\tan^2\eta)\,
\cos^p\xi\, F(-\frac {p-q}2, 
-\frac {p+q+6} 2, 4; -\tan^2\xi),  
\end{split}
\end{equation*}
for $x\in S$; see \cite{Johnson-2}.

\section {Restriction of $(SO_0(n, 1; \mathbb F), \pi_{\nu})$
to $(SO_0(n-1, 1; \mathbb F)$ }
\subsection{Principal series of $G$}

For $\nu\in \mathbb C$
let $\pi_{\nu}$ be the induced representation
of $G$ from $MAN$ consisting of measurable functions $f$ on $G$
such that
\begin{equation}
  \label{eq:ind-norm}
f(g me^{tH_0}
n
)=
e^{-\nu t}f(g),  me^{tH_0}n\in MAN  
\end{equation}
and $f\big{|}_{K}\in L^2(K)$. 
(Our  representation $\pi_\nu$ is $Ind_{MAN}^G(e^{\rho -\nu})$
in the standard notation \cite{Kn-book}. However
the parameter $\nu$ has some advantage it is ``stable''
under branching; see Theorem 3.6 below.) In particular $f$ in $\pi_\nu$
are invariant
under $M$, and $\pi_{\nu}$ is further  realized on 
$L^2(K/M)=L^2(S)$.
We denote $X_\nu$
 the corresponding $(\fg, K)$-module
with
$X_\nu$ the algebraic sum of $K$-irreducible subspaces
in $L^2(K/M)$. 
To simplify notation we
shall denote $\pi_{\nu}$ also
the corresponding unitary representation of $G$
when $(X_\nu,  \pi_{\nu}, \fg)$
is unitarizable.

The $L^2$-norm in $L^2(S)$
is not unitary
for $\pi_{\nu}$ except when $\nu=\rho+it$ for $t\in \mathbb R$,
$\rho$ being given by (\ref{eq:realrho}).
The unitarizable representations
 $(X_\nu,  \pi_{\nu}, \fg)$
for real $\nu$ are usually called complementary series.
They have been found 
in  \cite{Kostant-BAMS, Johnson-Wallach}. (See also
\cite{Dobrev-etal} for related results for the real group $SO_0(n,1)$.)
 The constant
$\lambda_{\nu}(\tau)$ below are rewritten in terms
of the Pochammer symbol $(a)_m$ and further the Gamma
functions.

\begin{theo+}
There is a positive definite
$\fg$-invariant form on $X_\nu$
given by
\begin{equation}
\label{nu-norm-gen}
\Vert w\Vert_{\nu}^2=\sum_{\tau} {\lambda_\nu(\tau)}
\Vert w_\tau\Vert^2, \quad w=\sum_{\tau} w_\tau \in 
X_\nu,
\end{equation}
where $\Vert w_\tau\Vert^2 $ is the $L^2$-norm, 
and its completion
forms an unitary irreducible representation
of $G$, if
 \begin{enumerate}
\item 
$\mathbb F=\br$, $0 <\nu < n-1$,
\begin{equation}
  \label{eq:lam-r}
\lambda_\nu(p)=\frac{(n-1-\nu)_p}
{(\nu)_p}
=\frac{\Gamma(n-1-\nu +p)}
{\Gamma(n-1-\nu) \Gamma(\nu +p)};  
\end{equation}
\item  $\mathbb F=\bc$,  $0 <\nu < 2n$,
\begin{equation}
  \label{eq:lam-c}
\lambda_\nu(p, q)=
\frac
{
(n-\frac \nu 2)_{p}
}
{(\frac \nu 2)_{p}
}
\frac
{
(n-\frac \nu 2)_{q}
}
{(\frac \nu 2)_{q}
}
=  
\frac
{\Gamma^2(\frac \nu 2)
\Gamma(n-\frac \nu 2 +p)}
{
\Gamma^2(n-\frac \nu 2 )  \Gamma(\frac \nu 2 +p)}
\frac
{\Gamma(n-\frac \nu 2 +q)}
{\Gamma(\frac \nu 2 +q)}
;
\end{equation}
\item  $\mathbb F=\bh$, $2 <\nu < 4n$,
\begin{equation}
  \label{eq:lam-h}
\begin{split}
\lambda_\nu(p, q)&=
\frac
{
(2n-\frac \nu 2)_{\frac{p-q}2}
}
{(\frac \nu 2 -1)_{\frac{p-q}2}
}
\frac
{
(2n+1-\frac \nu 2)_{\frac{p+q}2}
}
{(\frac \nu 2)_{\frac{p+q}2}
}\\
&
=
\frac
{\Gamma(\frac \nu 2 -1)\Gamma (\frac \nu 2)}
{\Gamma(2n-\frac \nu 2) \Gamma(2n+1-\frac \nu 2)}
\frac{
\Gamma
(2n-\frac \nu 2+\frac{p-q}2)
\Gamma(2n+1-\frac \nu 2+\frac{p+q}2)}
{\Gamma(\frac \nu 2 -1
+\frac{p-q}2)
\Gamma (\frac \nu 2 +\frac{p+q}2)
}.
\end{split}
\end{equation}
\end{enumerate}
\end{theo+}

\subsection{General criterion of boundedness}

We fix $n$ and let $H=SO_0(n-1, 1; \mathbb F) \subset G$
be the subgroup of elements of $g\in G$
fixing the  $n$-th coordinate $x_n$
in $\mathbb F^{n+1}$. Denote $L:=K\cap H$, a maximal subgroup
of $H$.
The subsphere, or the equator, $S^{d(n-1)-1}$ in $\mathbb F^{n-1}$
of the sphere $S=K/M \subset \mathbb F^{n}$
defined
by the equation $x_n=0$ will be
written as $S^\flat$, which is 
homogeneous space of $L$, $S^b=L/L\cap M$.
To avoid confusion
we denote by $\pi_\nu^\flat$ the 
corresponding  representations 
of $H$ and   $X_\nu^\flat$
the $L$-finite vectors ,
 and the corresponding
decomposition of $L^2(S^\flat)=L^2(L/L\cap M)$ will be
written as
$$
L^2(S^\flat)=\sum_{\sig }^\oplus
 V^{\sig}
$$
with $\sig$ being specified accordingly.

We shall need a general and elementary
criterion for boundedness
of intertwining operators.
The sufficient part of the following Lemma 3.2 is
used in \cite{Speh-Venk-1} implicitly, 
and we  give here a proof for the sake of completeness.
Let $K$  temporarily be a compact group
and  $L\subset K$ a closed
subgroup.
Let $(\mathcal W, \Vert \cdot \Vert_{\mathcal W})
$ and 
$(\mathcal V, \Vert \cdot \Vert_{\mathcal V})$ be 
unitary representations of  $K$ and respectively $L$.
Consider
$$
\mathcal W{\big |}_{K} =\sum_{\tau}^\oplus
 W^{\tau},\quad
\mathcal V{\big |}_{L} =\sum_{\sig}^\oplus V^{\sig}
$$
 the irreducible decomposition of $\mathcal W$
and $\mathcal V$
under $K$ and respectively $L$
counting multiplicities, all assumed  being finite. 
Consider further the branching of $\mathcal W^{\tau}$
under $L$. Write $\sig\subset \tau$ if a representation
$\sig$ appears in $\tau$ (counting multiplicities)
with $\widetilde V^{\tau, \sig}$
the corresponding irreducible component,
and denote $P_{\tau, \sig}
$ the corresponding
orthogonal projection, i.e., 
\begin{equation}
  \label{eq:-brchn-w-v}
 W^{\tau}
=\sum_{\sig\subset \tau}^\oplus
  \widetilde V^{\tau, \sig}, \quad 
P_{\tau, \sig}:  W^{\tau}\to \widetilde V^{\tau, \sig}.  
\end{equation}
Suppose $R$ is  a 
densely defined $L$-invariant operator
from $K$-finite elements in $\mathcal W$
to $L$-finite elements in $\mathcal V$, and
$$R_{\tau, \sig}:=P_{\tau, \sig} R: W^{\tau}\to 
V^{\sig}
$$
its  components, i.e.,  $R=\sum_{\tau}
\sum_{\sig\subset \tau} R_{\tau, \sig}$
on $K$-finite functions. We write $\Vert R\Vert_{\mathcal W, \mathcal V}$
its norm whenever it is finite.

\begin{lemm+}
\label{crit}
The restriction operator
$R$ extends to  a bounded operator
from $\mathcal W$ to 
$\mathcal V$ if and only if 
there is a constant $C$
such that for any $\sig$
\begin{equation}
  \label{eq:comp-n}
  \sum_{\tau\supset \sig}
\Vert R_{\tau, \sig}\Vert^2_{\mathcal W, \mathcal V}
 \le { C}.
\end{equation}
\end{lemm+}

\begin{proof} Let $w=\sum_{\tau} w_\tau
\in \mathcal W$ be an element with
finite many nonzero components $w_\tau$.
Its  squared norm in $\mathcal W$
is
 $$
\Vert w\Vert_{\mathcal W}^2
= \sum_{\tau} \Vert w_\tau\Vert_{\mathcal W}^2
$$
by our assumption.
 We compute the norm $\Vert Rw\Vert_{\mathcal V}$.
Writing $w$ as $w=\sum_{\tau} w_\tau
=\sum_{\tau} \sum_{\sig\subset \tau} 
P_{\tau, \sig} w_\tau
$,
we have
$$
Rw=\sum_{\sig} \sum_{\tau\supset \sig} 
R_{\tau, \sig}P_{\tau, \sig}
 w_\tau,
$$
and
$$
\Vert Rw\Vert_{\mathcal V}^2
=\sum_{\sig} 
\Vert \sum_{\tau\supset \tau} 
R_{\tau, \sig}P_{\tau, \sig}
 w_\tau\Vert_{\mathcal V}^2 
\le \sum_{\sig}
 (\Vert \sum_{\tau\supset \sig} 
\Vert R_{\tau, \sig}\Vert_{\mathcal W, \mathcal V}
 \Vert P_{\tau, \sig}
 w_\tau\Vert_{\mathcal W}
)
^2.
$$
If the condition (\ref{eq:comp-n}) is satisfied
we find, by Cauchy-Schwarz inequality, that
$$
\Vert Rw\Vert_{\mathcal V}^2\le \sum_{\sig} 
\left(
\sum_{\tau\supset \sig} 
\Vert R_{\tau, \sig}\Vert^2_{\mathcal W, \mathcal V}
\right) 
\left(\sum_{\tau\supset \sig} 
 \Vert P_{\tau, \sig}
 w_\tau\Vert^2_{\mathcal W}
\right) 
$$ 
which is dominated by
$$
C \sum_{\sig} \sum_{\tau\supset \sig} 
 \Vert P_{\tau, \sig}
 w_\tau\Vert^2_{\mathcal W}
=C \sum_{\tau}
\sum_{\sig\subset \tau} 
 \Vert P_{\tau, \sig}
 w_\tau\Vert^2_{\mathcal W}
=C \sum_{\tau}
\Vert w_\tau\Vert^2_{\mathcal W}
=C \Vert w\Vert_{\mathcal W}^2,
$$
finishing the proof of sufficiency. Conversely
suppose $R$ is a bounded operator. Then so is $R^\ast$,
and for a given $ v\in V^\sig$ we have
$$
C\Vert v\Vert_{\mathcal V}^2\ge
\Vert R^\ast v\Vert_{\mathcal W}^2
=\Vert \sum_{\tau\supset \sig}R^\ast_{\sig, \tau} v\Vert_{\mathcal W}^2
=\sum_{\tau\supset \sig}
\Vert 
R^\ast_{\sig, \tau} v
\Vert_{\mathcal W}^2.
$$
But each $R^\ast_{\sig, \tau}$ is a scalar constant
of an isometric operator by Schur's lemma,
 and we have 
$$\Vert 
R^\ast_{\sig, \tau} v
\Vert_{\mathcal W}^2=
\Vert 
R^\ast_{\sig, \tau}\Vert_{\mathcal V, \mathcal W}^2
\Vert v
\Vert_{\mathcal V}^2=
\Vert 
R_{\tau, \sig}\Vert_{\mathcal W, \mathcal V 
}^2
\Vert v
\Vert_{\mathcal V}^2.
$$
Substituting this into the above inequality
we obtain (\ref{eq:comp-n}).
\end{proof}

\subsection{Restriction of spherical harmonics}
We specify the above
considerations to the
restriction  $R: C^\infty(S)\to C^\infty(S^\flat), f(x', x_n)
\to f(x')$.
The branching of $W^{\tau}=\sum_{\sig}\widetilde
V^{\tau, \sig}
$ of an irreducible $K$-component $W^\tau$ under $L$
can  be read off  abstractly
from known results.
However
we need to find all isotypic $L$-irreducible
subspaces $\widetilde V^{\tau, \sig}\subset W^{\tau}$ 
 with nonzero restriction, i.e. with
the restriction
$$R_{\tau, \sig}: \widetilde V^{\tau, \sig}\to V^\sig
$$ acting as an isomorphism. 
More precisely we shall 
study the abstract branching    (\ref{eq:-brchn-w-v})
along with the concrete restriction
$$
W^\tau\big{|}_{x_n=0}:=\{g(x')=f(x', 0),\, x'\in
S^\flat; f\in W^\tau\}=\sum_{\sig\subset \tau} V^{\sig}.
$$
We shall
drop the upper-index $\tau$ in $\widetilde V^{\tau, \sig}
$ in
 the lemma below, as it is fixed in the summation.
The parameterization of $(\tau; \sigma)$ will be 
$(p; s)$ for $\mathbb F=\mathbb R$
and $(p, q; s, t)$ for $\mathbb F=\mathbb R, \mathbb C$.
Recall also the notation  in Lemma 2.1.
\begin{lemm+}
\begin{enumerate}
\item $\mathbb F=\mathbb R$. 
The branching of $W^p$ under 
 $L=SO(n-1)$  is multiplicity free. 
The restriction  $W^p\big{|}_{x_n=0}
$
under $L=SO(n-1)$ 
is decomposed as 
$$
W^p\big{|}_{x_n=0} =
\sum_{0\le s\le p, p-s \,\,\text{even}}^{\oplus}
V^{s}
$$
The corresponding unique $s$-isotypic 
component in $W^p$
is
given by (as functions on $S$)
$$
\tilde 
 V^{s}=
\{h(x') 
\phi_{p-s}^{n+2s}({x_n});
h\in V^{s}\}
$$
\item  $\mathbb F=\mathbb C$. The branching of $W^{p, q}$ under 
$L=U(n-1)$   is multiplicity free. 
The space  $W^{p, q}\big{|}_{x_n=0}$
under $L$ 
is decomposed as 
$$
W^{p, q}\big{|}_{x_n=0}=
\sum_{s\le p, t\le q, p-s=q-t}^\oplus
  V^{s, t}.
$$
For each $(s, t)$ the 
 unique $(s, t)$-isotypic 
component in $W^{p, q}$
is given by
$$
 \widetilde V^{s, t}=
\{h(x') 
\phi_{p-s, q-t}^{n+s+t}({x_n});
h\in V^{s, t}\}.
$$

\item  $\mathbb F=\mathbb H$.
The space  $W^{p, q}\big{|}_{x_n=0}$
under $L=Sp(n-1)\times Sp(1)$ 
is decomposed as 
$$
W^{p,q}\big{|}_{x_n=0} =
\sum_{0\le p-s \,\,\text{even} , t=q }^\oplus 
V^{s, t}.
$$
The corresponding  $(s, t)$-isotypic 
component is given by
$$
 \widetilde V^{s, t}=
\{h(z') \phi_{p-s, 0}^{n+\frac{s}2}({x_n});
h\in V^{s, t}\}
$$
\end{enumerate}
\end{lemm+}
\begin{proof} 
Let $\mathbb F=\mathbb R$.
The multiplicity
free result in this case  is well-known.
The statement on the restriction is
a result of Vilenkin \cite[(9), p.~495]{Vilenkin}.
The proof there relies on explicit
computations for the projection into spherical
harmonics, which seem not easy
to generalize to other cases. We give
a slightly different  proof 
which applies also to the other cases and which
avoids some redundant computations.
Denote $L_n=\sum_{j=1}^n
\frac{\partial^2}
{\partial x_j^2}$  the Laplacian on $\mathbb R^n$.
Recall that the spherical polynomial
 $f=r^m C_{m}^{\frac{n-2}2}(\frac{x_n}r)$ is
the unique $SO(n-1)$ invariant
polynomials on $\mathbb R^n$ of degree $m$ satisfying $L_n f=0$,
where $C_{m}^{\frac{n-2}2}(t)$ is the Gegenbauer polynomial.
Let $x=(x', x_n)\in \mathbb R^n$, and put $ u:=|x'|$,  $v:=x_n$.
We have 
$$
L_n =L_{n-1}+\frac{\partial^2}{\partial v^2}
=\frac{\partial^2}{\partial u^2}+
\frac{n-2}{u} \frac{\partial}{\partial u}+
\frac{\partial^2}{\partial v^2},
$$
when acting on functions depending only on $|x'|$ and $x_n$.
Rephrasing in terms of $u, v$
we have 
 the unique polynomial solution of the form
$f(u, v)=(u^2+v^2)^{\frac{m}2}
 C(\frac{v}{\sqrt{u^2+v^2}})
$,
of the equation 
\begin{equation}
  \label{eq:reduc-}
L_{n} f=\frac{\partial^2 f}
{\partial u^2} +
  \frac{n-2}{u}
\frac{\partial f}
{\partial u}  
 +
\frac{\partial^2 f}
{\partial v^2} =0.
\end{equation}
is when $C=C_{m}^{\frac{n-2}2}$.
Now for fixed $s\le p$ we
search for an isotypic $SO(n-1)$-component in 
$W^p$ of type $V^{s}$ consisting of homogeneous
polynomials $F(x)$ of degree $p$ of
the form $F(x)=h(x') f(u, v)
=h(x') f(|x|', x_n)$, where $h$
is a spherical harmonics of degree ${s}$ on $\mathbb R^{n-1}$, i.e.
$L_{n-1}h=0$, and
$f(u, v)=(u^2 +v^2)^{\frac{p-s}2}
C(\frac{v}{\sqrt{u^2 +v^2}})$.
The Laplace equation $L_n F=
(L_{n-1} +\frac{\partial^2}{\partial v^2})F=0$
becomes
$$
(L_{n-1}h(x'))  f(u, v)
+2\sum_{j=1}^{n-1}x_j\frac{\partial h(x')}
{\partial x_j}
\frac 1{u}
\frac{\partial h}
{\partial u}f(u, x_n)
+ h(x') L_n f(u, x_n)=0,
$$
with $L_n f$ computed in  (\ref{eq:reduc-}).
But $L_{n-1}h(x')=0$ and 
$\sum_{j=1}^{n-1}x_j\frac{\partial h(x')}
{\partial x_j}=s h(x')$ by our assumption. 
Thus it reduces
to 
\begin{equation}
2s\frac 1{u}
\frac{\partial h}
{\partial u}f(u, v)
+ L_n f(u, v)=0,
\end{equation}
i.e.
\begin{equation}
  \label{eq:reduc-1}
 \frac{\partial^2 f}
{\partial u^2} +
  \frac{n-2+2s}{u}
\frac{\partial f}{\partial u}  
+
\frac{\partial^2 f}
{\partial v^2} 
=0.
\end{equation}
This is precisely the equation
 (\ref{eq:reduc-}) with $n$ replaced by $n+2s$ and $m$
replaced by $p-s$. Thus
$f$ is a constant multiple of
$r^{p-s} C_{p-s}^{\frac{n-2}2+s}(\frac{x_n}r)$
(which is a posterior  polynomial in $x$). Note that
this is non-zero for $x_n=0$ only is $p-s$ is even. This
proves the case for $\mathbb F=\mathbb R$.

$\mathbb F=\mathbb C$. The multiplicity
free result is also known; see
e.g.  \cite{Kraemer-76}.
The abstract decomposition
of $
W^{p, q}\big{|}_{z_n=0}$
follows easily by counting the degrees $(s, t)$. 
We search an $L$-isotypic component
consisting of the polynomials of the form $F(x)=h(x')r^{p+q-s-t}C(\frac {x_n}r)$
as above. This
results  to characterization
of  the function $C(\frac {x_n}r)$
exactly the same as $\phi_{\tau}^m$
as in Lemm 2.1 (2) with  $\tau=(p', q')$
and $m$ determined by $(p, q)$ and $n$.

$\mathbb F=\mathbb H$. The group $Sp(1)$ acts on the space
of polynomials on the right, $h\in Sp(1):f(x)\mapsto f(xh)$,
and it acts on the  space $W^{p, q}$  as the symmetric
tensor $\odot^q(\mathbb C^2)$.
So does it
also on the space $W^{p, q}\big{|}_{x_n=0}$. Thus
any irreducible component
must be of type $V^{s, t}$ with  $t=q$,
 again by (2.1).
In particular $p-s=(p-q)-(s-t)$ is even since
both $p-q$ and $s-t$ are even.
This proves the decomposition. 
The rest of the proof is almost
the same as above. (Note that the function
$\phi_{p-s, 0}^{n+\frac{s}2}$ is obtained
from $\phi_{p, 0}^{n}$ in Lemma 2.1 (3)
by formally replacing $n$ by
$n+\frac{s}2$, which is not necessarily an integer.)
\end{proof}

We compute now the norm of $R_{\tau, \sig}$.
For positiva constants  $C_{\tau, \sig}$
and $D_{\tau, \sig}$ we write 
$C_{\tau, \sig}\sim D_{\tau, \sig}$ if both 
$\frac{C_{\tau, \sig}}{D_{\tau, \sig}}$
and $\frac{D_{\tau, \sig}}{C_{\tau, \sig}}
$ are dominated by positive constants independent of $\tau, \sig$.

\begin{prop+} With the notation as above we 
have the $L^2(S)-L^2(S^\flat)$-norm of $R_{\tau, \sig}: W^\tau\to V^\sig
\subset L^2(S^\flat)$ is given by
\begin{enumerate}
\item $\mathbb F=\mathbb R$, $p-s\ge 0 $ even,
$$
\Vert R_{p,  s}\Vert^2= 
\frac{\Gamma(\frac n 2)}
{\Gamma(\frac {n-1} 2)\Gamma(\frac 1 2)
} 
\frac
{
(2p+n-2)
\Gamma(\frac{n+p+s-2}2)
\Gamma(\frac{p-s+1}2)
}
{
\Gamma(\frac{p-s+2}2)
\Gamma(\frac{n+p+s-1}2)
}
\sim
 \frac{p+1}
{(p+s+1)^{\frac 12} 
(p-s+1)^{\frac 12}
};
$$
\item  $\mathbb F=\mathbb C$, $p\ge s\ge 0$, 
$q\ge t\ge 0$, $p-q=s-t$,
$$
\Vert R_{(p, q),  (s, t)}\Vert^2=   p+q+n-1;
$$
\item  $\mathbb F=\mathbb H$, 
 $p-s=2k\ge 0$ and $s-t\ge $ even, $q=t$,
\begin{equation*}
  \begin{split}
\Vert R_{(p,q), (s, t)}\Vert^2&=  \frac{\Gamma(2n-2)}
{\Gamma(2n)}
(k+1)(2k+2(n-1)+s -1)(k+2(n-1)+s)\\
&\sim (k+1)(k+s +1)^2.    
  \end{split}
\end{equation*}
\end{enumerate}
In all other cases of $(\tau, \sig)$ we have $R_{\tau, \sig}=0$.
\end{prop+}
\begin{proof} $\mathbb F=\mathbb R$. 
 By the previous lemma
and Schur lemma
we see that  $R_{p, s}: W^p\to V^{s}$ is up to a constant
a partial isometry, and  $\widetilde{V^{s}}\to V^{s}$ is up to a constant
an isometry. Thus 
$$
\Vert R_{p, s}\Vert^2=\frac{\Vert Rf\Vert^2}{\Vert f\Vert^2}
$$
for any $0\ne f\in \widetilde{V^{s}}$. 
Now we take $f=h(x')
\phi_{p-s}^{n+2s}(x_n)
$, which has a form of variable separation,
 and we have, by  (\ref{eq:polar})
and the $s$-homogeneity of $h(x')$, that
$$
\Vert f\Vert^2 =
\frac{\Gamma(\frac n 2)}
{\Gamma(\frac {n-1} 2)
\Gamma(\frac 1 2)} 
\int_{|x_n|< 1} 
(1-|x_n|^2)^{\frac {n-3}2 +s}
|\phi_{p-s}^{n+2s}(x_n)|^2
\int_{S^\flat}  |h(y')|^2 dy'
dx_n,
$$
and 
\begin{equation}
  \label{eq:Rf}
\Vert Rf\Vert^2
=|\phi_{p-s}^{n+2s}(0)|^2
 \int_{S^\flat}  |h(y')|^2 dy'.  
\end{equation}
Consequently
$$
\Vert R_{p, s}\Vert^2=
|\phi_{p-s}^{n+2s}(0)|^2 
\left(\frac{\Gamma(\frac n 2)}
{\Gamma(\frac {n-1} 2)
\Gamma(\frac 1 2)} 
\int_{|x_n|< 1} 
(1-|x_n|^2)^{\frac {n-3}2 +s}
|\phi_{p-s}^{n+2s}(x_n)|^2 dx_n \right)^{-1}.
$$
Note that the integral 
$$I:=\int_{|x_n|< 1}
(1-|x_n|^2)^{\frac{n-3}2 +s}
|\phi_{p-s}^{n+2s}(x_n)|^2 dx_n 
$$ 
is up to a constant the square norm in $L^2(S^{n+2s-1})$
of the spherical polynomial
$\phi_{p-s}^{n+2s}(x_n)$ in dimension $n+2s$, and can be evaluated
by using
  (\ref{eq:schur-1}) in terms
of the dimension $\text{dim} W^{p-s}_{n+2s}
$
of the representation of $SO(n+2s)$.
 The
exact (a rather subtle)  constant is
computed in   (\ref{eq:polar}),
$$
\int_{|x_n|< 1} 
(1-|x_n|^2)^{\frac{n-3}2 +s}
|\phi_{p-s}^{n+2s}(x_n)|^2dx_n 
= \frac{\Gamma(\frac{n+2s-1}2)\Gamma(\frac 12)}
{\Gamma(\frac {n+2s}2)}
        \frac 1 {\text{dim} W^{p-s}_{n+2s}}.
$$
Thus using  the dimension formula that
$$
\text{dim}W_n^{j}=\binom{n+j-1}{j}-\binom{n+j-3}{j-2}
=
\frac
{
(n+2j-2)\Gamma(n+j-2)
}
{\Gamma(j+1) \Gamma(n-1),
}
$$
we have $I^{-1}$ is 
\begin{equation*}
I^{-1} =\frac{\Gamma(\frac {n+2s}2)}
{\Gamma(\frac{n+2s-1}2)\Gamma(\frac 12)}
\frac{(2p+n-2)\Gamma(n+p+s-2) 
 }
{\Gamma (n+2s-1) \Gamma(p-s+1)}.  
\end{equation*}

The evaluation $\phi_{p-s}^{n+2s}(0)$ 
 in (\ref{eq:Rf}) is zero unless
 $p-s=2k$ is even, in which case it is
$$
(-1)^{k}
\frac{
(-k)_{k}
(-\frac{p-s-1}2)_{k}
}
{(\frac{n+2s-1}2)_{k} 
k!}.
$$
But $(-k)_k=(-1)^k k!$, 
$(-\frac{p-s-1}2)_{k}=(-1)^k
(\frac 12)_{k}
=(-1)^k \frac{\Gamma(\frac 12 +k) }{\Gamma(\frac 12)}
$, we find that the evaluation,
disregarding the sign $(-1)^k$
and the constant $\Gamma(\frac 12)$, is
$$
\frac
{\Gamma(\frac{p-s+1}2) \Gamma(\frac{n+2s-1}2)
}
{\Gamma(\frac{n+p+s-1}2)
}.
$$
Using the product formula $\Gamma(2x) =\Gamma(\frac 12)^{-1}2^{2x-1}
{\Gamma(x)}{\Gamma(x+\frac 12)}$ 
we obtain then the
formula for $\Vert R_{p, s}\Vert^2$
as stated. The rest follows
from the Stirling formula that
$$
\frac{\Gamma(n+a)}{\Gamma(n+b)}\sim n^{a-b}, \quad n\to \infty.
$$

The case $\mathbb F=\bc$ 
is done by similar
computations. 
In the case $\mathbb F=\bh$ 
we have
$$
\Vert R_{(p, q), (s, q)}\Vert^2=
\frac{
|\phi_{p-s}^{n+\frac{s}2}(0)
|^2 \int_{S^\flat} |h(x')|^2 dx'
}
{\Vert \phi_{p-s}^{n+\frac{s}2} h\Vert^2
}
$$
with ${\Vert \phi_{p-s}^{n+\frac{s}2} h\Vert^2
}$ being
$$
\frac{2\Gamma(n)}{\Gamma(
n-1)
\omega_{3}}
\int_{x_n\in \bh, |x_n|<1}
|\phi_{p-s}^{n+\frac{s}2}(x_n)|^2
(1-|x_n|^2)^{\frac 12( 4(n-1)-2 +  2s)}
\int_{S^\flat} |h(x')|^2 dx' dx_n
$$
by the integral formula above
for separated variables.
The norm of $\phi^{n+\frac {s}2}$ can not
be computed using the dimension formula
for $s$ odd as it can not be interpreted
as spherical polynomials on a symmetric
space. However we may use
by known integral formulas \cite{Askey-etal-bk,Szeg}
for Jacobi polynomials
$P^{(\alpha, \beta)}(t)$
on the interval $[-1, 1]$.  (More generally
one may use the theory
of Heckman-Opdam \cite{Heckman-Opdam-1} for Jacobi polynomials
with general root multiplicities.)
Indeed the function $\phi_{k, 0}^n$ in \S2.2
for any real $n> 1$ can be written
as
$$
\phi_{k, 0}^n(x)=
\frac{\Gamma(k+1)\Gamma(2n-2)}{
\Gamma(k+2n-2)}
P_{k}^{(2n-3, 1)}(2|x|^2-1)
$$
where $|x|$ is the norm of a quaternionic number $x\in \bh$.
The norm to be computed is
$$
\int_{x\in \bh, |x|<1}
|\phi_{m, 0}^m(x)|^2(1-|x|^2)^{\frac 12(4(n-1)-2} dx
=\omega_3 
\int_{0}^1
|\phi_{m, 0}^m(x)|^2(1-|x|^2)^{\frac 12(4(n-1)-2} 
x^3 dx
$$
and which is further \cite[(6.4.5)-(6.4.6), pp. 299-301]{Askey-etal-bk}
(see also \cite{Szeg})
$$
\omega_3 \frac{\Gamma^2(2n-2)\Gamma(k+1)\Gamma(k+2)
}
{
\Gamma(k+2n-2)\Gamma(k+2n-1) (2k+2n-1)
}.
$$
The rest is done by a routine computation.
\end{proof}

Note that when $\mathbb F=\br$ and $n=3$ our result coincides with
that in \cite[Lemma 2.4]{Speh-Venk-1}. For $\mathbb F=\bc$,
and $W^{p, q}=W^{p, 0}$ the space
the holomorphic polynomials of degree $p$, the norm
of $R$ can be found directly by computing
of the integral $\int_S |x_1^p|^2dx$ on the  sphere $S$
in $\mathbb C^n$.

\subsection{Discrete components of complementary series}
Before stating our first main result 
we note the following
elementary 
\begin{lemm+}
Suppose $0<\alpha<1, \beta >0$, $\alpha +\beta >1$
and $\gamma>1$. Then
\begin{equation*}
  \label{eq:el-su}
\sum_{j=0}^\infty \frac 1{(j+1)^\alpha (q+j+1)^\beta}
\le C \frac{1}{q^{\alpha+\beta -1}},   \quad 
\sum_{j=0}^\infty \frac {1}{ (j+q+1)^\gamma}
\le C \frac{1}
{
(q+1)^{\gamma-1}
}, \quad \forall q\ge 0
\end{equation*}
\end{lemm+}
The second estimate is straightforward.
The first sum is dominated by the integral
$$\int_{0}^\infty 
 \frac{1}{x^\alpha (x+q+1)^{\beta}}dx
=\frac{1}{(q+1)^{\alpha+\beta-1}}
\int_{0}^\infty 
\frac{1}{x^\alpha (x+1)^{\beta}}dx
=\frac{1}{(q+1)^{\alpha+\beta-1}} C
$$
since the integral  $\int_{0}^\infty \frac{1}{x^\alpha (x+1)^{\beta}}dx = C<\infty$
is convergent by our assumption.

Observe also that 
$$
R: (X_\nu, \pi_{\nu}, \fg)\to 
(X^\flat_\nu, \pi_{\nu}^\flat, \fh), \,
f(x)\mapsto f(x', 0)
$$
intertwines the action of $\pi_{\nu}^\flat$ of  $\fh$. 
Thus
the boundedness of $R$ implies that  $(\pi_{\nu}^\flat, \fh)$
is a discrete component whenever both are unitarizable. 
In accordance with the notation
$\Vert \cdot \Vert$ in Theorem 3.1
we denote $\Vert T\Vert_{\nu, \mu}$ the norm
of an operator $T: X_\nu \to X_\mu^\flat$
and further $\Vert T\Vert_{\nu}=\Vert T\Vert_{\nu, \nu}$.
We have then
$$
\Vert R_{\tau, \sig}
\Vert_{\nu, \mu}^2
=\frac{\lambda_\mu(\sig)^\flat}{\lambda_\nu(\tau)} \Vert R_{\tau, \sig}\Vert^2,
$$
and the criterion (\ref{eq:comp-n}) becomes
\begin{equation}
  \label{eq:comp-n-pre}
  \sum_{\tau\supset \sig}\Vert R_{\tau, \sig}\Vert^2
\lambda_\nu(\tau)^{-1}
 \le \frac{ C}{   \lambda_\mu^\flat(\sig)}.
\end{equation}

\begin{theo+} The restriction of $(\pi_{\nu}, G)$
on $H$ contains $(\pi_{\nu}^\flat, H)$
as a discrete component in the following cases
\begin{enumerate}
\item $\mathbb F=\mathbb R$, $n\ge 3$, $0<\nu<\frac{n-2}2$;
\item $\mathbb F=\mathbb C$, $n\ge 3$, $0<\nu<n-2$;
\item $\mathbb F=\mathbb H$, $n\ge 2$, $2<\nu<2n-1$.
\end{enumerate}
\end{theo+} 
\begin{proof} 
 $\mathbb F=\mathbb R$. 
First note that
$\frac{n-2}2< n-2 <n-1$ thus both
$(\pi_{\nu}, G)$ and $(\pi_{\nu}, H)$
are well-defined unitary \reps/.
 We use  now Lemma 3.2 with $\tau=p$ and $\sig=s$.
 The constants $\lambda_{\nu}(p)$,  $\lambda^\flat_{\mu}(s)$ 
and the series (\ref{eq:comp-n-pre}) in question are
$$
\lambda_\nu(p)\sim {(p+1)^{n-1-2\nu}},
\quad \lambda_\nu^\flat (s)\sim {(s+1)^{n-2-2\nu}},
$$
$$
\sum_{p\ge s, p-s\, \text{even}}
 \frac{p+1}
{(p+s+1)^{\frac 12} 
(p-s+1)^{\frac 12}
}\frac 1 {(p+1)^{n-1-2\nu}}.
$$
Writing $p=s+2j$ we see the sum is dominated
by
\begin{equation*}
\begin{split}
\sum_{j=0}^\infty
 \frac{s+2j+1}
{(2s+2j+1)^{\frac 12} 
(2j+1)^{\frac 12}
}
\frac 1 {(s+2j+1)^{n-1-2\nu}}
\le 
C \sum_{j=1}^\infty
 \frac{1}
{j^{\frac 12}
}
\frac 1 {
(s+j)^{n-1-2\nu-\frac 12}},
\end{split}
\end{equation*}
and further
 by $(s+1)^{-(n-2-2\nu)}$
in view of Lemma 3.5, namely by $\frac 1{\lambda^\flat(s)}$.

$\mathbb F=\mathbb C$. $\lambda_{\nu}$ has the asymptotics
$$
\lambda_\nu(p, q)\sim (p+1)^{n-\nu} (q+1)^{n-\nu}.
$$
For fixed type $(s, t)$ of $L$
 the series $ \sum_{\tau\supset \sig}
\Vert R_{\tau, \sig}\Vert^2
\lambda(\tau)^{-1}$
is  dominated up to a constant by
\begin{equation*}
  \begin{split}
&\quad \sum_{p-s=q-t\ge 0}
\frac {p+q+2}
{(p+1)^{n-\nu} (q+1)^{n-\nu} }\\
&
= \sum_{p-s=q-t\ge 0}
(\frac 1{(p+1)^{n-\nu-1} (q+1)^{n-\nu} }
+\frac 1{(p+1)^{n-\nu} (q+1)^{n-\nu-1} }  )
  \end{split}
\end{equation*}
as sum of two, say $I+II$. Now
$$
I=\sum_{k=0}^\infty \frac 1{
(s +k+1)^{n-\nu -1}
(t +k+1)^{n-\nu}
},
$$
and 
$$
I\le \frac 1{(s +1)^{n-\nu-1}}
\sum_{k=0}^\infty \frac 1{
(t +k+1)^{n-\nu}
}\le C \frac 1{(s +1)^{n-\nu-1}
(t +1)^{n-\nu-1}
}\le C\frac 1{\lambda_\nu^\flat (s, t)}
$$
by Lemma 3.5. The same holds for $II$.

$\mathbb F=\mathbb H$. Writing $p=s+2k$, $k\ge 0$,
we have
$$
\lambda_\nu(p, q)\sim 
({p-q}+1)^{2n+1-\nu} 
(p+q+1)^{2n+1-\nu}
\sim ({s-q}+k+1)^{2n+1-\nu} 
(s+q+k+1)^{2n+1-\nu}
$$
and
$$\Vert R_{(p,q), (s, q)}\Vert^2  \sim
(k+1)(s+k +1)^2.
$$
The sum  (\ref{eq:comp-n-pre}) is bounded by
\begin{equation*}
\begin{split}
&\quad\,\sum_{k=0}^\infty
\frac{(k+1)(k+s +1)^2}
{(s-q+k+1)^{2n+1-\nu}
(s+q+k+1)^{2n+1-\nu}
}
\\
&
\le
\sum_{k=0}^\infty
\frac{k+1}
{(s-q+k+1)^{2n+1-\nu}
(s+q+k+1)^{2n-1-\nu}
}
\\
&
\le
\frac 1{(s+q+1)^{2n-1-\nu}}
\sum_{k=0}^\infty
\frac{k+1}
{
(s-q+k+1)^{2n+1-\nu}
}
\\
&\le
\frac 1{(s+q+1)^{2n-1-\nu}}
\sum_{k=0}^\infty
\frac{1}
{
(s-q+k+1)^{2n-\nu}
}
\\
&\le  C
\frac 1{(s+q+1)^{2n-1-\nu}}
\frac{1}
{(s-q+1)^{2n-\nu-1}
}\sim \frac 1{\lambda_{\nu}^\flat({s, t})},
\end{split}
\end{equation*}
finishing the proof.
\end{proof}

\begin{rema+}
For  $\mathbb F=\mathbb R$ and $n=3$
the full decomposition of the complementary series
$\pi_\nu$ of $SO_0(3,1)$ under $SO_0(2, 1)$
is done in \cite{Mukunda}. If (in terms of
our parametrization) $\frac 12\le \nu <1$
it is a sum of two direct integrals of spherical principal
series, and if $0<\nu <\frac 12$ there is one extra
discrete component,  the complementary series.
 The appearance
of complementary series $\pi_{\nu}^\flat$
of  $SO_0(n-1, 1)$ in the complementary series
$\pi_{\nu}$ of $SO_0(n, 1)$ is done in \cite{Speh-Venk-2}
using the non-compact realization on $\mathbb R^{n-1}$.
It might be possible to find
a full decomposition for general
$G=SO_0(n, 1;\mathbb F)$ using the techniques
in \cite{Mukunda}. The most interesting
part might still be the 
the discrete spectrum in view 
of its stability under restriction
and induction in the  Ramanujan duals \cite{Burger-Sarnak, Burger-Li-Sarnak}.
\end{rema+}

\subsection{The quotients
 $(\mathcal W, \pi_{\nu})$
at negative integers $\nu$
and their discrete components}

The representation $ \pi_{\nu}$ are reducible 
\cite{Johnson-Wallach}
for $\nu$ satisfying
certain integral conditions, and there exist
 unitarizable subrepresentations
(or quotients). More precisely we have
the following result \cite{Takahashi-mathfr, Johnson-Wallach},
retaining the notation of $\lambda_\tau$
as the Schur proportional constants; 
here we have rewritten 
them in similar formulation as in Theorem 3.1.

\begin{theo+}
There is a unitarizable
 irreducible quotient 
$(\mathcal W_{\nu}, \pi_{\nu})$, 
whose completion forms  a unitary irreducible 
representation of  $G$ in the following cases
\begin{enumerate}
\item $\mathbb F=\br$, $n\ge 3$, $\nu= -k$, 
$k\ge 0$, $\mathcal W_{\nu}
=X_{\nu}/M_{\nu}$,
$$
 M_\nu=\sum_{p=0}^k W^p,
$$
$$
\lambda_{\nu}(p)=
\frac{(n-1-\nu +k+1)_{p-k-1}}
{(\nu +k+1)_{p-k-1}}
=C_\nu
\frac{\Gamma(n-1-\nu +p)}
{\Gamma(\nu +p)}
;
$$
\item $\mathbb F=\bc$,  $n\ge 2$,
$\nu=-2k$, $k> 0$,
 $\mathcal W_{\nu}
=X_{\nu}/M_{\nu}$,
$$
M_\nu=
\sum_{p\le k, q\ge 0} W^{p,q} +
\sum_{q\le k, p\ge 0} W^{p,q}
\quad (k>0),
$$
\begin{equation*}
  \begin{split}
\lambda_\nu(p, q)
&=
\frac{(n-\frac \nu 2 +k+1)_{p-k-1}
(n-\frac \nu 2 +k+1)_{q-k-1}}
{
(\frac \nu 2 +k+1)_{p-k-1}
(\frac \nu 2 +k+1)_{q-k-1}}
\\
&=C_\nu
\frac
{
\Gamma(n-\frac \nu 2 +p)
\Gamma(n-\frac \nu 2 +q)}
{
\Gamma(\frac \nu 2 +p)
\Gamma(\frac \nu 2 +q)
}
  \end{split}
\end{equation*}
and for $k=0$ with three quotients 
$(\mathcal W_{0}^\pm, \pi_0^\pm)$,
$(\mathcal W_{0}, \pi_0)$,
$$\mathcal W_{0}^+
=\sum_{p=0 }^\infty W^{p, 0}/
\mathbb C, \quad
\mathcal W_{0}^-
=\sum_{q=0}^\infty W^{0, q}/
\mathbb C
$$
$$\lambda_{\nu}^+(p)=\frac{\Gamma(p)}{\Gamma(n+p)},
\lambda_{\nu}^-(q)=\frac{\Gamma(q)}{\Gamma(n+q)},
$$
and 
$$
\mathcal W_{0}=X_0/\sum_{p=0}^\infty(W^{p, 0}+
W^{0, p})
$$
$$
\lambda_0(p, q)=\frac{\Gamma(n+p-1)\Gamma(n+q-1)}
{\Gamma(p)\Gamma(q)};
$$

\item 
 $\mathbb F=\bh$, $n\ge 1$, $\nu=-2k$, $k \ge -1$,
$\mathcal W_{\nu}
=X_{\nu}/M_\nu$, 
$$
 M_\nu=
\sum_{p-q\le 2k+2} W^{p,q}, \quad k\ge 0, \quad
 \quad M_{\nu}=W^{0, 0}, \quad k= -1,
$$
\begin{equation*}
\begin{split}
\lambda_\nu(p, q)&=
\frac
{
(2n-\frac \nu 2 +k+1)_{\frac{p-q}2 -k-1}
(2n+1-\frac \nu 2 +k+1)_{\frac{p+q}2 -k-1}
}
{
(\frac \nu 2-1 +k+1)_{\frac{p-q}2 -k-1}
(\frac \nu 2 +k+1)_{\frac{p+q}2 -k-1}
}\\
&
=C_\nu \frac
{\Gamma(2n-\frac \nu 2 +\frac{p-q}2)
\Gamma(2n+1-\frac \nu 2 +\frac{p+q}2)}
{
\Gamma(\frac \nu 2-1 +\frac{p-q}2)
\Gamma(\frac \nu 2 +\frac{p+q}2)
}.
\end{split}
\end{equation*}
\end{enumerate}
\end{theo+}

Note that the same $\nu=\nu(k)$
as above is also a reducible point
for $(X_{\nu, \flat}, \pi_{\nu}^\flat, \fh)$.
The corresponding
quotient \rep/ for $\fh$
will be written as
 $(\mathcal V_{k}, \pi_{\nu(k)}^\flat, \fh)$

\begin{theo+} Let $n\ge 4$ for $\mathbb F=\br$,
$n\ge 3$ for $\mathbb F=\bc$,
and $n\ge 2$ for $\mathbb F=\bh$.
 The representation
$(\mathcal V_{k}, 
\pi_{\nu(k)}^\flat, \fh)$
(and the corresponding completion as  for $H$)
appears as an irreducible discrete component
in  $(\mathcal W_{k}, \pi_{\nu(k)}, \fg)$ 
(respectively  for $G$)
 restricted to $\fh$ (resp. $H$).
\end{theo+}
\begin{proof}
Let $Q=Q_{\nu}$ 
be the quotient map
$
Q: X_{\nu}^{\flat}\to
X_{\nu}^{\flat}/M_{\nu}^\flat:=\mathcal V_{k}$
at  the reducible point
$\nu$ as above for the group $H$. The map
$
QR:  X_\nu \to X_{\nu}^{\flat}\to 
\mathcal V_{k}
$
is clearly  $(\pi_{\nu}^\flat, \fh)$ intertwining
 and induces
a map
$$
QR:  \mathcal W^k= X_{\nu}/M_{\nu}
 \to X_{\nu}^{\flat}/M_{\nu}^\flat=\mathcal V_{k}.
$$
We prove the boundedness of $QR$ by the method
above. Notice that the asymptotics exponent
of $\lambda(\nu)$ has the same
dependence for positive $\nu$, e.g.
in the  case $\mathbb F=\br$ with $\nu=-k$,
$$
\lambda_\nu(p)\sim  {(p+1)^{n-1-2\nu}}, \quad p\ge k+1,
$$
and $n-1-2\nu \ge 2$. Thus the same 
proof carries over to all cases, and we
omit the details.
\end{proof} 

There is some slightly difference
when $n=3$ for $\mathbb F=\br$,
as $H=SO_0(2,1)$
 has its maximal compact subgroup
being the torus and there is a splitting
of the restriction to holomorphic and 
antiholomorphic discrete series. Note
that we have also excluded the case
 $\mathbb F=\bc$, $n=2$, namely $SU(2,1)$,
as the  restriction map above  is zero
on the quotient $\mathcal W_{-2k}$; actually
$\mathcal W_{-2k}$ is 
a discrete series and 
it branching under $SU(1, 1)$
can possibly be studied using some
general tools \cite{Ors-Var, kob-rest-2005}.

\subsection{The representation $\pi_0$
and $\pi_0^\pm$ for $SU(n, 1)$}
The representation $\pi_{}^\pm$ 
on the quotient
$$
\pi_0^+ =\sum_{p=0}^\infty W^{p, 0}/\mathbb C,
\quad
\pi_0^- =\sum_{p=0}^\infty W^{0, p}/\mathbb C
$$
are unitarizable representation of $\fg$. 
 $\pi_0^{+}$ can be constructed  also
by using  the analytic continuation 
of the weighted Bergman space \cite[Theorem 5.4]{FK}, i.e. scalar
holomorphic discrete series, on the unit ball 
$G/K$ in $\mathbb C^n$
with reproducing kernel
$(1-(z, w))^{-\mu}$ at the reducible point $\mu=0$;
see e.g. \cite{Hw-Liu-Z}
where a reproducing kernel and its expansion are found
for the space.
A full decomposition under $SU(n-1, 1)$ of the series and their quotient
can be obtained easily. Indeed let
 $\pi_0^{\pm, \flat}$ be the corresponding
\rep/ for $H$ and $\pi_{j}^{+, \flat}$
the unitary \rep/ of $H$
realized as the space of holomorphic
functions on the unit ball $\{z\in \bc^{n-1}, |z|<1\}$
with
reproducing kernel $(1-(z, w))^{j}$,  with  $H$
acting as
$$
g=\begin{bmatrix} a&b\\
c&d
\end{bmatrix}\in H, f(z)\mapsto (cz+d)^{-j}f((az+b)(cz+d)^{-1}).
$$
$\pi^{+}_{j, \flat}$ 
is a discrete series of $H$
only when $j\ge n$. Define analogously 
$\pi^{-, \flat}_{j}$ in terms of conjugate holomorphic
functions.
The following result follows from easy consideration of
expansion of holomorphic functions $f(x)$ in 
the last variable $x_n$.

\begin{prop+}
The \rep/ $(\pi_0^{\pm}, G)$ is decomposed under $H$
as
$$
\pi_0^\pm =\pi_0^{\pm, \flat} \oplus 
(\sum_{j=1}^\infty{}^\oplus
\pi^{\pm, \flat}_{j}).
$$
\end{prop+}

\section {Restriction of $(F_{4(-20)}, \pi_{\nu})$
to $H=Spin(8, 1)$ }

\subsection{The subgroup $Spin(8, 1)$}
Recall that  $H_0\in \fp=\mathbb O^2$
has nonzero roots $\pm 2, \pm 1$ in $\fg$. Denote
$\fg_{\pm 2}$ and $\fg_{\pm 1}$
the respective root spaces. Then the Lie algebras
$\fg_{\pm 1}$ generate a subalgebra of
$\fg$ of rank one which is easily seen
to be $\fh:=spin(8,1)$. The Cartan
decomposition of $\fh$ is
$\fh=spin(8)\oplus\mathbb O
$ with 
$spin(8)$ acting on $\mathbb O$ by the spin representation.
The simply connected
subgroup of $G$ with Lie algebra $\fh$ is then
$H=Spin(8, 1)$ whose maximal compact group
is $L=Spin(8)$; see e.g. \cite{Burger-Li-Sarnak}.
 It follows from the decomposition (\ref{eq:r16-m})
that the stabilizer of $H_0\in \mathbb O\subset \fh$
in $H$ is  also $M=Spin(7)$ and that 
$L/M=Spin(8)/Spin(7)$ is the  sphere $S^7$.
We have thus
$$
L^2(S^7)=\sum_{p\ge 0}^\oplus
 V^p
$$
where $V^p$ is the space of spherical harmonics
of degree $p$ on $S^7$, defined by the condition $x_2=0$
in $S^{15}=\{x=(x_1, x_2)\in \mathbb O^2; |x|=1\}$.
The  decomposition of $L^2(S^{15})=L^2(K/M)$
is given in (\ref{eq:L-2-cayley})
with $(W^{p, q})^M=\mathbb C\phi_{p, q}$.
We consider now 
the restriction of $W^{p, q}\big{|}_{x_2=0}$
of the components $W^{p, q}$.

\begin{lemm+} 
The decomposition of $W^{p, q}$
under $Spin(8)$ is multiplicity free and 
 $W^{p, q}\big{|}_{x_2=0}=V^{q}$; in other words
the only irreducible component
in the decomposition with non-zero restriction
to $S^7$ is the representation $V^{q}$. Moreover
the square norm of $R: W^{p, q}\to V^q$ is given
by
\begin{equation*}
  \begin{split}
\Vert R\Vert^2
&=C\frac{(p+7)\prod_{j=0}^2(p+q+8+2j)(p-q+2+2j)(q+4+2j)(q+1+2j)}
{(q+3)(q+1)_5}
 \\
&\sim {(p+1) (p+q+1)^3(p-q+1)^3}
  \end{split}
\end{equation*}
where $C$ is a numerical constant independent of $p$ and $q$.
\end{lemm+}

\begin{rema+} The \rep/  $
W^{p, q}
$ of $K=Spin(9)$
is of highest weight
 $\frac{p-q}2\lambda_1 + q\lambda_4$
with $\lambda_1, \lambda_2, \lambda_3,  
\lambda_4$
 the fundamental weights dual
to the simple roots $\alpha_1,
\alpha_2, \alpha_3, \alpha_4$. In the standard
notation
they are  $\alpha_1=e_1-e_2$,
$\alpha_2=e_2-e_3$,
$\alpha_3=e_3-e_4$,
$\alpha_4=e_4$ and $\frac{p-q}2\lambda_1 + q\lambda_4
=\frac p2e_1 +\frac q2 e_2 +\frac q2 e_3 +\frac q2e_4
=\frac 12(p, q, q, q)
$. 
The four simple roots for $spin(8)$ are
$\delta_1=e_1-e_2, \delta_2=e_2-e_3, \delta_3=e_3-e_4, \delta_4
=e_3+e_4$. 
The above branching rule above can be formulated
as 
$$
W^{\frac 12(p, q, q, q)}\big{|}_{x_2=0}
=V^{\frac q2(1, 1, 1, 1)
},
$$
with $V^{\frac q2(1, 1, 1, 1)}$
being the space of spherical harmonics of degree $q$ on
$S^7$, which is also   of highest weight $q e_1$
as $SO(8)$ representation.
This discrepancy of highest weights is explained by the triality
in $Spin(8)$.
The Dynkin diagram of $Spin(8)$
is
  \begin{center}
\begin{figure}[h]
\includegraphics{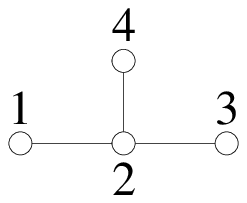}
\end{figure}
  \end{center}
There is a symmetry of $S_3$ (as outer automorphisms)
acting on the  three simple roots $\delta_1, \delta_3, \delta_4$.
The highest weight $\frac q2(1, 1, 1, 1)
=\frac q2(\delta_1+2\delta_2+\delta_3+2\delta_4)$,
whereas $qe_1=\frac q2( 2\delta_1+2\delta_2+\delta_3+\delta_4 )$
and the permutation $(134)$ exchanges the two weights.
Also the multiplicity one property of
$W^{p, q}
$ under $M=Spin(7)$ factors through
$Spin(8)$ and we have
$(W^{p, q})^M{\big |}_{x_2=0}
=(W^{p, q}{\big |}_{x_2=0})^M
=(V^{q})^M
=\phi^{8}_q$. Note that $M=Spin(7)$ in $Spin(8)\subset Spin(9)$
is not the obvious copy of $Spin(7)$ in $Spin(9)$
defined by the standard inclusion $\mathbb R^7\subset
\mathbb R^8\subset \br^9$; 
in the space $W^{p, q}$ the former
copy $M=Spin(7)$ has only one-dimensional fixed vectors,
whereas the
latter copy $Spin(7)$ has
 arbitrarily large multiplicities by
the construction of Gelfand-Zetlin basis \cite{Vilenkin}.
\end{rema+} 

\begin{proof}The first statement is well-known.
Any irreducible representation
of $Spin(8)$ in $W^{p, q}\big{|}_{x_2=0}$
is a constituent in $L^2(S^7)$
and contains thus a unique
$M=Spin(7)$-invariant element. But
$(W^{p, q}\big{|}_{x_2=0})^M=
(W^{p, q})^M\big{|}_{x_2=0} 
=\bc \phi^{p, q}\big{|}_{x_2=0}$,
and $\phi^{p, q}\big{|}_{x_2=0}$ is
$$
\phi_{p, q}
(\cos\eta, 0)=\cos^q \eta
F(-\frac q 2, -\frac {q-1} 2; \frac 72; -\tan^2 \eta)
$$
which is precisely the $M$-invariant spherical
harmonics  $
\phi^{8}_q(\cos 
\eta)$ on  $S^7$,
Section 2.2. Thus $W^{p, q}\big{|}_{x_2=0}$ is nonzero
and is just $V^{q}$. In particular the element
$\phi_{p, q}$ is in the isotypic component $\widetilde
V^{q}\subset W^{p, q}$ of $V^q$. The squared norm of $R$ on $W^{p, q}$ is
$$
\Vert R\Vert^2=\Vert R \phi_{p, q}\Vert^2
\Vert^2 \phi_{p, q}
\Vert^{-2}, \quad 
R \phi_{p, q}=\phi^{8}_q.
$$
Both norms can be evaluated by the dimension formula.
Following the notation in the above remark
we have $W^{p, q}$ has highest weight
$\frac p2e_1 +\frac q2 e_2 +\frac q2 e_3 +\frac q2e_4$
with the positive roots being $\{e_i \pm e_j, e_i, \, 1\le i<j\le 4\}
$, and the dimension of
$W^{p, q}$ is then
\begin{equation*}
  \begin{split}
\text{dim} W^{p, q}&=
C_1 (p+7) \prod_{j=0}^2 (p+q+8+2j)(p-q+2+2j)(q+4+2j)(q+1+2j)\\
&
\sim (p+1)(p+q+1)^3(p-q+1)^3 (q+1)^6
 \end{split}
\end{equation*}
whereas the dimension of $V^{q}$ is
$$
\text{dim} V^{q}
=C_2 (q+3)(q+1)_5\sim (q+1)^6
$$
for some  constants $C_1, C_2$ independent of $p$ and $q$.
This completes the proof.
\end{proof}

\subsection {Discrete components}

Define the principal series
representation
$\pi_\nu$ of $G$ as 
in (3.1), realized on $L^2(K/M)=L^2(S^{15})$.
We recall  results in \cite{Johnson-2}
 on the  spherical complementary
series of $G$.

\begin{theo+} 
Let $6 <\nu <16$.
There is a positive definite
$(\fg, \pi_\nu)$-invariant form 
on the $(\fg, K)$-module 
$\sum_{p, q}W^{p,q}$
defined by
as in (\ref{nu-norm-gen}) with
$$
\Vert w_{p, q}\Vert_{\nu}^2={\lambda_\nu(p, q)}\Vert w_{p, q}\Vert^2,
\quad \lambda_\nu(p, q)=\frac
{
(8-\frac \nu 2)_{\frac{p-q}2}
}
{(\frac \nu 2- 3)_{\frac{p-q}2}
}
\frac
{
(11-\frac \nu 2)_{\frac{p+q}2}
}
{
(\frac \nu 2)_{\frac{p+q}2}
}.
$$
\end{theo+}

We study now the branching of the complementary
series under $H=Spin(8, 1)$.
Denote $\pi_{\mu}^\flat$
the principal series representation
of $SO_0(8, 1)$, thus also for $H$,  as defined for 
in   (\ref{eq:ind-norm}).
Note that there is
a discrepancy between the normalization of $H_0$
here and there;
the roots of $H_0$
in $\fh$ is $\pm 2$ here instead of $\pm 1$ as
in Section 2.1. 
In particular,
the restriction map $R: f(x_1, x_2)\mapsto Rf(x_1)=f(x_1, 0)$, $(x_1, x_2)\in \mathbb O^2$,
defines an $\fh$-intertwining operator
$$
R: (X_{\nu}, \pi_\nu, \fg)
\to (X_{\frac\nu 2}^{\flat}, \pi_{\frac \nu 2}^\flat, \fh).
$$

\begin{theo+} Let $6 <\nu <7$.
The restriction of $(\pi_{\nu}, G)$
on $H$ contains $(\pi_{\frac \nu 2}^\flat, H)$
as a  discrete component.
\end{theo+}

\begin{proof} The $\lambda_\nu$ and $\lambda_{\mu}^\flat$
in this case are
$$
\lambda_\nu(p, q)\sim (p-q+1)^{11-\nu}(p+q+1)^{11-\nu},
\quad \lambda_{\mu}(q)\sim  (q+1)^{7-2\mu}
$$
with $p-q=2k\ge 0$ even. The sum to be treated is
  \begin{equation*}
 \begin{split}
&\quad\,\sum_{k=0}^\infty \frac{q+k+1}
{(k+1)^{8-\nu}
(2k+q+1)^{8-\nu}
}
\\
&
\le \sum_{k=0}^\infty \frac{1}
{(k+1)^{8-\nu}
(2k+q+1)^{7-\nu}
}
\\
&
\le
\frac 1{(q+1)^{7-\nu}}
\sum_{k=0}^\infty \frac{1}
{(k+1)^{8-\nu}
}
\\
&=
C\frac 1{(q+1)^{7-\nu}}
\sim \frac 1{ \lambda_{\frac \nu 2}^\flat (q)},
 \end{split}
  \end{equation*}
completing the proof.
\end{proof}

\begin{rema+} The complementary series in 
is parametrized in 
 \cite[Example C]{Burger-Li-Sarnak} as $-5\le \lambda\le  5$, i.e. the standard
parametrization \cite{Kn-book}. Our $\nu$ is their
$\rho +\lambda=11+\lambda$. It is stated
there that the point $\lambda=3$, i.e. $\nu=8$ is
in the automorphic dual $\hat G_{\text{aut}}$ of $G$. Note that
this point falls outside the range $6 <\nu<7$ in our
theorem. One can deduce from the 
Burger-Li-Sarnak conjecture on the Ramanujan
dual $\hat H_{\text{Raman}}$
 for $H=SO(n, 1)$
and our theorem above
on some nonexistence of certain intervals in the 
set $\hat H_{\text{Raman}}$.
In view of  \cite[Theorem 1]{Burger-Li-Sarnak}
it would be also interesting to study 
the induction of automorphic representations
of $H$ to $G$.
\end{rema+}

The representation  $\pi_{\nu}$
has also unitarizable subquotients at integral $\nu$:
for $\nu=6-2k$, $k\ge 0$, the quotient
$$
\mathcal W^{\nu}=X_\nu/{M_\nu}, \quad
M_\nu =\sum_{p-q\le k}W^{p, q},
$$
is unitarizable; see \cite{Johnson-2}.
However in this case 
the restriction composed
with quotient map  is zero.
Presumable there is no discrete
component
under $H$ and it would be  interesting
to pursue this further.

The main results in the present paper prove 
the existence of one single discrete component
under $H$ of  a complementary series
of $G$. It may happen that there
are more  discrete components.
 We shall study them in a forthcoming paper.
Applying  the stability result
 \cite[Theorem 1 (ii)]{Burger-Li-Sarnak} 
to  those cases where $\pi_{\nu}^\flat$ appear discretely in $\pi_{\nu}$
we conclude finally
if  $\pi_{\nu}$ is an automorphic representation
of  $G$
so is  $\pi_{\nu}^\flat$.


\def\cprime{$'$} \newcommand{\noopsort}[1]{} \newcommand{\printfirst}[2]{#1}
  \newcommand{\singleletter}[1]{#1} \newcommand{\switchargs}[2]{#2#1}
  \def\cprime{$'$} \def\cprime{$'$} \def\cprime{$'$}
\providecommand{\bysame}{\leavevmode\hbox to3em{\hrulefill}\thinspace}
\providecommand{\href}[2]{#2}

\end{document}